# An analytic solution to the Busemann-Petty problem on sections of convex bodies

By R. J. Gardner, A. Koldobsky, and T. Schlumprecht*

**Abstract**

We derive a formula connecting the derivatives of parallel section functions of an origin-symmetric star body in $\mathbb{R}^n$ with the Fourier transform of powers of the radial function of the body. A parallel section function (or $(n-1)$-dimensional X-ray) gives the ($(n-1)$-dimensional) volumes of all hyperplane sections of the body orthogonal to a given direction. This formula provides a new characterization of intersection bodies in $\mathbb{R}^n$ and leads to a unified analytic solution to the Busemann-Petty problem: Suppose that $K$ and $L$ are two origin-symmetric convex bodies in $\mathbb{R}^n$ such that the ($(n-1)$-dimensional) volume of each central hyperplane section of $K$ is smaller than the volume of the corresponding section of $L$; is the ($n$-dimensional) volume of $K$ smaller than the volume of $L$? In conjunction with earlier established connections between the Busemann-Petty problem, intersection bodies, and positive definite distributions, our formula shows that the answer to the problem depends on the behavior of the $(n-2)$-nd derivative of the parallel section functions. The affirmative answer to the Busemann-Petty problem for $n \leq 4$ and the negative answer for $n \geq 5$ now follow from the fact that convexity controls the second derivatives, but does not control the derivatives of higher orders.

## 1. Introduction

The 1956 Busemann-Petty problem (see [BP]) asks the following question. Suppose that $K$ and $L$ are origin-symmetric convex bodies in $\mathbb{R}^n$ such that

$$\mathrm{vol}_{n-1}(K \cap H) \leq \mathrm{vol}_{n-1}(L \cap H)$$

---

*This work was partly carried out when the second author participated in the Workshop in Linear Analysis and Probability Theory at Texas A&M University, 1997. First author supported in part by NSF Grant DMS-9501289; second author supported in part by NSF Grant DMS-9531594 and the UTSA Faculty Research Award; third author supported in part by NSF Grants DMS-9501243 and DMS-9706828 and by the Texas Advanced Research Program Grant 160766.

1991 *Mathematics Subject Classification*. Primary: 52A20; secondary: 46B07, 42B10.

for every hyperplane $H$ containing the origin; does it follow that

$$\text{vol}_n(K) \leq \text{vol}_n(L)?$$

The concept of an intersection body of a star body was introduced by Lutwak [Lu] in 1988 and played a crucial role in the solution to the Busemann-Petty problem. A slightly more general notion was defined in [GLW], as follows. An origin-symmetric star body $K$ in $\mathbb{R}^n$ is said to be an *intersection body* if there exists a finite (non-negative) Borel measure $\mu$ on the $(n-1)$-dimensional sphere $S^{n-1}$ so that the radial function $\rho_K$ of $K$ equals the spherical Radon transform of $\mu$ (all necessary definitions will be given in §3). If the measure $\mu$ has a continuous positive density on $S^{n-1}$, then there exists another star body $L$ so that the radial function of $K$ at every point $\xi \in S^{n-1}$ is equal to the $(n-1)$-dimensional volume of the section of $L$ by the hyperplane $\xi^\perp = \{x \in \mathbb{R}^n : \langle x, \xi \rangle = 0\}$ (see [Ga3, Ch. 8] for details). In this case $K$ is said to be the *intersection body of a star body*.

A connection between the Busemann-Petty problem and intersection bodies was established by Lutwak [Lu], and slightly modified in [Ga2] and [Z1], [Z2]. In particular, [Z1, Th. 2.22] is as follows.

THEOREM A. *The Busemann-Petty problem has an affirmative answer in $\mathbb{R}^n$ if and only if every origin-symmetric convex body in $\mathbb{R}^n$ is an intersection body.*

However, it turned out to be quite difficult to calculate the inverse spherical Radon transform of a radial function in order to check if a given body is an intersection body.

The Busemann-Petty problem has a long history. A negative answer to the problem for $n \geq 5$ was established in a sequence of papers by Larman and Rogers [LR] (for $n \geq 12$), Ball [Ba] ($n \geq 10$), Giannopoulos [Gi] and Bourgain [Bo] ($n \geq 7$), Papadimitrakis [P], Gardner [Ga1] and Zhang [Z2, Th. 6.1] ($n \geq 5$). A little later, it was proved in [Ga2] that every origin-symmetric convex body in $\mathbb{R}^3$ is an intersection body and, therefore, that the answer to the Busemann-Petty problem is affirmative when $n = 3$ (note that the answer is trivially affirmative when $n = 2$). The result followed from the fact that, when $n = 3$, the inverse spherical Radon transform $R^{-1}\rho_K$ of the radial function of an origin-symmetric strictly convex body with $C^\infty$ boundary satisfies

$$R^{-1}\rho_K(\xi) = -\frac{1}{4\pi^2} \int_0^\infty \frac{A'_\xi(z)}{z} dz,$$

for all $\xi \in S^{n-1}$, where

$$A_\xi(z) = \text{vol}_{n-1}(K \cap \{x \in \mathbb{R}^n : \langle x, \xi \rangle = z\})$$

is the parallel section function of $K$ in the direction $\xi$.

In 1997, it was shown by the second named author [K4] that an origin-symmetric cube in $\mathbb{R}^4$ is an intersection body. The result in [K4] was a consequence of the following connection between intersection bodies and the Fourier transform, established in [K3].

THEOREM B. *An origin-symmetric star body $K$ in $\mathbb{R}^n$ is an intersection body if and only if the radial function $\rho_K$ is a positive definite distribution on $\mathbb{R}^n$.*

Theorem B has other applications. For example, it was shown in [K3] that the unit ball of every finite-dimensional subspace of $L_p$, $0 < p \leq 2$, is an intersection body, which, in particular, confirms the conjecture of Meyer [M] that the answer to the Busemann-Petty problem is affirmative if the body $K$ is a polar projection body (unit ball of a subspace of $L_1$). The paper [K5] presents a variety of examples of origin-symmetric convex bodies in $\mathbb{R}^n$, $n \geq 5$, that are not intersection bodies. Theorem B provides an effective method for determining whether a given star body is an intersection body; in fact it follows from [K2, Lemma 1], that

$$2^n \pi^{n-1} R^{-1} \rho_K = \widehat{\rho_K}$$

when the Fourier transform $\widehat{\rho_K}$ of $\rho_K$ is a continuous function on $S^{n-1}$.

After learning the results of [K3], [K4], Zhang [Z3] proved that every origin-symmetric convex body in $\mathbb{R}^4$ is an intersection body, which implies an affirmative answer to the Busemann-Petty problem when $n = 4$. The proof is based on a geometric argument, similar to that of [Ga2], which shows that if $K$ is an origin-symmetric convex body with $C^2$ boundary, then

$$R^{-1} \rho_K(\xi) = -\frac{1}{16\pi^2} A''_\xi(0),$$

for all $\xi \in S^{n-1}$. It is then an immediate consequence of the Brunn-Minkowski theorem that the inverse spherical Radon transform is non-negative.

In this article, we establish the following formula.

THEOREM 1. *Let $K$ be an origin-symmetric star body in $\mathbb{R}^n$ with $C^\infty$ boundary, and let $k \in \mathbb{N} \cup \{0\}$, $k \neq n - 1$. Suppose that $\xi \in S^{n-1}$, and let $A_\xi$ be the corresponding parallel section function of $K$.*

(a) *If $k$ is even, then*

$$(\rho_K^{n-k-1})^\wedge(\xi) = (-1)^{k/2} \pi(n - k - 1) A_\xi^{(k)}(0);$$

(b) *if $k$ is odd, then*

$$(\rho_K^{n-k-1})^\wedge(\xi) = c_k \int_0^\infty \frac{A_\xi(z) - A_\xi(0) - A''_\xi(0)\frac{z^2}{2} - \cdots - A_\xi^{(k-1)}(0)\frac{z^{k-1}}{(k-1)!}}{z^{k+1}} \, dz,$$

where $c_k = (-1)^{(k+1)/2} 2(n-1-k)k!$, $A_\xi^{(k)}(0)$ is the derivative of order $k$ of the function $z \mapsto A_\xi(z)$ at zero, and $(\rho_K^{n-k-1})^\wedge$ is the Fourier transform in the sense of distributions.

Note that in Theorem 1(b) all the derivatives of $A_\xi$ of odd order vanish, since $A_\xi(\cdot)$ is even. In particular, $A'_\xi(0) = 0$; using this, integration by parts, and the above equation relating $R^{-1}\rho_K$ and $\widehat{\rho_K}$, we see that the formula of [Ga2] given above is just the case $n = 3$, $k = 1$ of Theorem 1(b). Furthermore, the formula of [Z3] is the special case $n = 4$, $k = 2$ of Theorem 1(a). The theorem therefore represents a generalization of these earlier formulas to arbitrary dimensions.

We apply Theorem 1(a) to confirm that the answer to the Busemann-Petty problem is affirmative when $n = 4$, and use Theorem 1(b) with $n = 5$ and $k = 3$ to present a simple example that confirms the negative answer when $n \geq 5$ (see §2). Therefore, Theorem 1, in conjunction with Theorems A and B, provides a unified analytic solution to the Busemann-Petty problem. Moreover, Theorems 1 and B give a characterization of intersection bodies in higher dimensions. For example, if $n$ is even, then an origin-symmetric star body $K$ in $\mathbb{R}^n$ with $C^\infty$ boundary is an intersection body if and only if the $(n-2)$-nd derivative of the function $(-1)^{(n-2)/2} A_\xi$ at zero is non-negative for every $\xi \in S^{n-1}$. Note that Theorem 1 with $k = 0$ gives the Fourier transform formula for the volume of central hyperplane sections of $K$, which was used in [K2] to confirm the conjecture of Meyer and Pajor on the minimal sections of $\ell_p^n$-balls with $0 < p \leq 2$. Putting $k = n$ and using the fact that an origin-symmetric convex body $K$ is a zonoid if and only if $\|x\|$ is a negative definite function (see [Le, pp. 219–223]), one gets a new characterization of zonoids.

The proof of Theorem 1 will be given in §4. We shall first use the concept of fractional derivatives to extend the mapping

$$k \mapsto \left.\frac{\partial^k}{\partial t^k} A_\xi(t)\right|_{t=0}, \quad k \in \mathbb{N},$$

to an analytic function

$$q \mapsto A_\xi^{(q)}(0), \quad \text{where } q \in \mathbb{C}, \quad \operatorname{Re} q > -1, \quad q \neq n-1,$$

and show that this extension satisfies the following formula.

THEOREM 2. *Let $K$ be an origin-symmetric star body in $\mathbb{R}^n$ with $C^\infty$ boundary and Minkowski functional $\|\cdot\|$. Suppose that $\xi \in S^{n-1}$, and let $A_\xi$ be the corresponding parallel section function of $K$. For $q \in \mathbb{C}$ with $\operatorname{Re} q > -1$, $q \neq n-1$,*

$$A_\xi^{(q)}(0) = \frac{\cos \frac{q\pi}{2}}{\pi(n-q-1)} (\|x\|^{-n+q+1})^\wedge(\xi).$$

The construction of $A_\xi^{(q)}$ together with necessary definitions and properties of distributions will be given in §3.

We conclude the introduction by formulating the isomorphic Busemann-Petty problem: Does there exist an absolute (not depending on the dimension) constant $c$ such that $\mathrm{vol}_n(K) \leq c\,\mathrm{vol}_n(L)$ whenever $\mathrm{vol}_{n-1}(K \cap H) \leq \mathrm{vol}_{n-1}(L \cap H)$ for every hyperplane $H$ containing the origin? This question is equivalent to the famous hyperplane (or slicing) problem, which remains one of the most important unsolved mysteries of the local theory of Banach spaces (see [MP]).

The results of this paper were announced in [GKS].

## 2. Applications of Theorem 1 to intersection bodies and the Busemann-Petty problem

We first prove that the answer to the Busemann-Petty problem is affirmative when $n \leq 4$. In view of Theorem A, it is enough to show that every origin-symmetric convex body in $\mathbb{R}^n$, $n \leq 4$ is an intersection body. Also, since the intersection of an intersection body with a hyperplane $H$ containing the origin is also an intersection body in $H$ (see [FGW], [GW], or [Z4, Lemma 3]; one can also deduce it from Theorem B), it is enough to consider the case $n = 4$. The following theorem was first proved by Zhang [Z3].

THEOREM 3. *Every origin-symmetric convex body $K$ in $\mathbb{R}^4$ is an intersection body.*

*Proof.* A result of Zhang [Z1, Th. 2.13] implies that an origin-symmetric convex body that is not an intersection body can be approximated arbitrarily closely in the Hausdorff metric by origin-symmetric convex bodies with $C^\infty$ boundaries that are also not intersection bodies. Therefore we can assume that $K$ has $C^\infty$ boundary. Put $n = 4$ and $k = 2$ in Theorem 1. We get $\widehat{\rho_K}(\xi) = -\pi A_\xi''(0)$ for every $\xi \in \mathbb{R}^4 \setminus \{0\}$. By the Brunn-Minkowski theorem (see, for example, [S, Th. 6.1.1]), the function $A_\xi$ is log concave, and, since $A_\xi(\cdot)$ is even, we see that $A_\xi''(0) \leq 0$ for every $\xi$. Thus, $\rho_K$ is a positive definite distribution. The result follows from Theorem B. $\square$

In view of Theorem A and the remark at the beginning of this section, the following theorem confirms the result of [P], [Ga1], and [Z2, Th. 6.1] that the answer to the Busemann-Petty problem is negative when $n \geq 5$. Note that the proofs in [Ga1] and [Z2, Th. 6.1] were based on the fact that certain special origin-symmetric bodies in $\mathbb{R}^5$ are not intersection bodies. The simple proof given here follows quickly from the case $n = 5$ and $k = 3$ of Theorem 1(b).

THEOREM 4. *There is an origin-symmetric convex body $K$ in $\mathbb{R}^5$ that is not an intersection body.*

*Proof.* Put $n = 5$ and $k = 3$ in Theorem 1. By Theorem B, it is enough to make sure that there exists a $\xi \in S^4$ such that the parallel section function $A_\xi$ of $K$ satisfies

$$\int_0^\infty \frac{1}{z^4} \left( A_\xi(z) - A_\xi(0) - A_\xi''(0)\frac{z^2}{2} \right) dz < 0.$$

To this end, let $\varepsilon \in (0, 1)$, define $f_\varepsilon(x) = (1 - x^2 - \varepsilon x^4)^{1/4}$, and let $a_\varepsilon > 0$ be such that $f_\varepsilon(a_\varepsilon) = 0$ and $1 - x^2 - \varepsilon x^4 > 0$ on $(0, a_\varepsilon)$.

The function $f_\varepsilon$ has its maximum at 0 and

$$f_\varepsilon''(x) = -(\frac{1}{2} + 3\varepsilon x^2)(1 - x^2 - \varepsilon x^4)^{-3/4} - 3(-\frac{1}{2}x - \varepsilon x^3)^2(1 - x^2 - \varepsilon x^4)^{-7/4},$$

so $f_\varepsilon$ is strictly concave on $[0, a_\varepsilon]$. It follows that

$$K = \left\{ (x_1, \ldots, x_5) \in \mathbb{R}^5 : x_5 \in [-a_\varepsilon, a_\varepsilon] \text{ and } \left( \sum_{i=1}^4 x_i^2 \right)^{1/2} \leq f_\varepsilon(|x_5|) \right\}$$

is a strictly convex body. Since for $0 \leq z \leq a_\varepsilon$,

$$K \cap \{(x_1, \ldots, x_5) \in \mathbb{R}^5 : x_5 = z\}$$

is a 4-dimensional Euclidean ball of radius $f_\varepsilon(z)$, we deduce that when $\xi = (0, 0, 0, 0, 1)$,

$$A_\xi(z) = \frac{\pi^2}{2} f_\varepsilon(z)^4 = \frac{\pi^2}{2}(1 - z^2 - \varepsilon z^4),$$

for $0 \leq z \leq a_\varepsilon$. This implies that the above integral equals $-\varepsilon a_\varepsilon \pi^2/2 < 0$. $\square$

Putting $k = n - 2$ in Theorem 1 and using Theorem B, we get a characterization of intersection bodies in $\mathbb{R}^n$ in terms of the derivatives of parallel section functions. In particular, if $n$ is even, a star body with $C^\infty$ boundary is an intersection body if and only if $(-1)^{(n-2)/2} A_\xi^{(n-2)}(0) \geq 0$ for every $\xi \in S^{n-1}$. This observation yields an informal explanation of the answer to the Busemann-Petty problem: Convexity implies that the parallel section functions are log concave (a property involving the first and second derivative), but does not provide any control over the third and higher derivatives.

## 3. Notation and auxiliary facts

The *spherical Radon transform* is the bounded linear operator on $C(S^{n-1})$ defined by

$$Rf(\xi) = \int_{S^{n-1} \cap \xi^\perp} f(x)\,dx, \qquad f \in C(S^{n-1}), \quad \xi \in S^{n-1}.$$

(Here and throughout, differentials such as $dx$ denote integration with respect to the Hausdorff measure of the appropriate dimension.) If $\mu$ is a finite Borel measure on $S^{n-1}$, then the spherical Radon transform of $\mu$ is defined as a measure $R\mu$ on $S^{n-1}$ such that, for every $f \in C(S^{n-1})$,

$$\langle R\mu, f \rangle = \langle \mu, Rf \rangle = \int_{S^{n-1}} Rf(\xi) \, d\mu(\xi).$$

Let $\phi$ be an integrable function on $\mathbb{R}^n$ also integrable on hyperplanes, let $\xi \in S^{n-1}$, and let $t \in \mathbb{R}$. Then

$$\mathcal{R}\phi(\xi; t) = \int_{\langle x, \xi \rangle = t} \phi(x) \, dx$$

is the *Radon transform of $\phi$ in the direction $\xi$ at the point $t$*. Now for arbitrary $\xi \in \mathbb{R}^n \setminus \{0\}$, the Radon transform in the direction of $\xi$ at $t$ is defined by

$$\mathcal{R}\phi(\xi; t) = \frac{1}{\|\xi\|_2} \mathcal{R}\phi\left(\frac{\xi}{\|\xi\|_2}; \frac{t}{\|\xi\|_2}\right),$$

where $\|\cdot\|_2$ is the Euclidean norm. By the well-known connection between the Fourier transform and the Radon transform (see [H, p. 4], where the notation is different), it follows that for every $\xi \in \mathbb{R}^n \setminus \{0\}$ and $s \in \mathbb{R}$,

(1) $$\hat{\phi}(s\xi) = \big(\mathcal{R}\phi(\xi; t)\big)^{\wedge}(s),$$

where on the right-hand side we have the Fourier transform of the function $t \mapsto \mathcal{R}\phi(\xi; t)$.

Let $K$ be a body that is star-shaped with respect to the origin. The *radial function* of $K$ is given by

$$\rho_K(x) = \max\{a > 0 : ax \in K\}, \quad x \in \mathbb{R}^n \setminus \{0\}.$$

We call $K$ a *star body* if $\rho_K$ is continuous *and positive* on $S^{n-1}$ (there are different definitions of this term in the literature; in particular, it is often not assumed that $K$ contains the origin in its interior).

Let $K$ be an origin-symmetric star body. We denote by $\|x\|_K = \min\{a > 0 : x \in aK\}$ the *Minkowski functional* on $\mathbb{R}^n$ generated by $K$. Clearly, $\rho_K(x) = \|x\|_K^{-1}$. In the sequel, $\|\cdot\| = \|\cdot\|_K$ will always denote the Minkowski norm of $K$.

For every $\xi \in S^{n-1}$, we define the *parallel section function* $z \mapsto A_\xi(z)$, $z \in \mathbb{R}$ of $K$ by

$$A_\xi(z) = \operatorname{vol}_{n-1}\big(K \cap (\xi^\perp + z\xi)\big) = \int_{\langle x, \xi \rangle = z} \chi(\|x\|) \, dx = \mathcal{R}\chi(\|x\|)(\xi; z),$$

where $\chi$ is the indicator function of $[-1, 1]$. (The function $A_\xi(z)$ is sometimes called the $(n-1)$-dimensional X-ray orthogonal to $\xi$; see [Ga3, Ch. 2].) For an arbitrary $\xi \in \mathbb{R}^n \setminus \{0\}$, we put

$$A_\xi(z) = \frac{1}{\|\xi\|_2} \mathcal{R}\chi(\|x\|) \left( \frac{\xi}{\|\xi\|_2}; \frac{z}{\|\xi\|_2} \right).$$

Then, by (1), for every fixed $\xi \in \mathbb{R}^n \setminus \{0\}$, the Fourier transform of the function $z \mapsto A_\xi(z)$ is equal to

(2) $$\widehat{A_\xi}(t) = \bigl(\chi(\|x\|)\bigr)^{\wedge}(t\xi).$$

The main tool of this paper is the Fourier transform of distributions. We use the notation from [GS]. As usual, we denote by $\mathcal{S}$ the space of rapidly decreasing infinitely differentiable functions on $\mathbb{R}^n$ with values in $\mathbb{C}$. By $\mathcal{S}'$ we signify the space of distributions over $\mathcal{S}$. The *Fourier transform* of a distribution $f$ is defined by $\langle \hat{f}, \hat{\varphi} \rangle = (2\pi)^n \langle f, \varphi \rangle$ for every test function $\varphi$. If a test function $\varphi$ is even,

$$(\hat{\varphi})^{\wedge} = (2\pi)^n \varphi \qquad \text{and} \qquad \langle \hat{f}, \varphi \rangle = \langle f, \hat{\varphi} \rangle$$

for every $f \in \mathcal{S}'$. If $q$ is not an integer, then the Fourier transform of the function $|z|^q$, $z \in \mathbb{R}$, is equal to (see [GS, p. 173])

(3) $$(|z|^q)^{\wedge}(t) = -2\Gamma(1+q)\sin\frac{q\pi}{2}|t|^{-q-1}, \quad t \in \mathbb{R}.$$

A distribution $f$ is called *positive definite* if, for every test function $\varphi$,

$$\left\langle f, \varphi * \overline{\varphi(-x)} \right\rangle \geq 0.$$

A distribution is positive definite if and only if its Fourier transform is a positive distribution (in the sense that $\langle f, \varphi \rangle \geq 0$ for every non-negative test function $\varphi$; see, for example, [GV, p. 152]).

For $t \in \mathbb{R}$, let $t_+ = \max\{0, t\}$. If $\varphi \in \mathcal{S}$ vanishes on a neighborhood of 0, the integral

$$\langle t_+^\lambda, \varphi(t) \rangle = \int_0^\infty t^\lambda \varphi(t)\, dt$$

exists for all $\lambda \in \mathbb{C}$, and, moreover, the function

$$\lambda \mapsto \int_0^\infty t^\lambda \varphi(t)\, dt$$

is complex differentiable on $\mathbb{C}$, and thus an analytic function. We now regularize the functional $\varphi \mapsto \langle t_+^\lambda, \varphi \rangle$, in order to define it on all of $\mathcal{S}$, in the following way (cf. [GS, Ch. I, § 3]). For $\lambda \in \mathbb{C}$ and $m \in \mathbb{N}$ such that $-m - 1 < \operatorname{Re}\lambda$, $\lambda \neq -1, -2, \ldots, -m$, and for every $\varphi \in \mathcal{S}$, we put

(4)
$$\langle t_+^\lambda, \varphi(t) \rangle = \int_0^1 t^\lambda \left( \varphi(t) - \varphi(0) - t\varphi'(0) - \cdots - \frac{t^{m-1}}{(m-1)!}\varphi^{(m-1)}(0) \right) dt$$
$$+ \int_1^\infty t^\lambda \varphi(t)\, dt + \sum_{k=1}^m \frac{\varphi^{(k-1)}(0)}{(k-1)!(\lambda + k)}.$$

If $-m-1 < \operatorname{Re}\lambda < -m$, we have

(5)
$$\langle t_+^\lambda, \varphi(t)\rangle = \int_0^\infty t^\lambda \left(\varphi(t) - \varphi(0) - t\varphi'(0) - \cdots - \frac{t^{m-1}}{(m-1)!}\varphi^{(m-1)}(0)\right) dt,$$

since in this case

$$\int_1^\infty t^{\lambda+k-1} dt = -\frac{1}{\lambda+k},$$

for $k = 1, \ldots, m$. The family $\{t_+^\lambda : \lambda \in \mathbb{C} \setminus \{-1, -2, \ldots\}\}$ forms an analytic distribution ([GS, p. 48]); that is, for any $\varphi \in \mathcal{S}$, the function $\lambda \mapsto \langle t_+^\lambda, \varphi(t)\rangle$ is analytic on $\lambda \in \mathbb{C} \setminus \{-1, -2, \ldots\}$. Furthermore, $\langle t_+^\lambda, \varphi(t)\rangle$ has, for each $k \in \mathbb{N}$, a simple pole at $\lambda = -k$ with residue $\varphi^{(k-1)}(0)/(k-1)!$ (see [GS, p. 49]). The function $\lambda \mapsto \Gamma(\lambda+1) = \int_0^\infty t^\lambda e^{-t}\, dt$ also has, for each $k \in \mathbb{N}$, a simple pole at $\lambda = -k$ with residue $(-1)^{k-1}/(k-1)!$. We conclude that

$$\left\{\frac{t_+^\lambda}{\Gamma(\lambda+1)} : \lambda \in \mathbb{C} \setminus \{-1, -2, \ldots\}\right\}$$

can be extended to an analytic distribution on $\mathbb{C}$, still denoted by $\{t_+^\lambda/\Gamma(\lambda+1) : \lambda \in \mathbb{C}\}$; and for $\lambda = -k$ and $\varphi \in \mathcal{S}$,

$$\left\langle \frac{t_+^\lambda}{\Gamma(\lambda+1)}, \varphi(t)\right\rangle = (-1)^{k-1}\varphi^{(k-1)}(0)$$

(see [GS, p. 56]). Outside any neighborhood of 0 the functional $t_+^\lambda/\Gamma(\lambda+1)$ acts like a finite measure, so that we can actually apply $t_+^\lambda/\Gamma(\lambda+1)$ to any continuous function that is infinitely differentiable on a neighborhood of 0, deducing the same conclusions as for functions in $\mathcal{S}$.

Such a function is the parallel section function $z \mapsto A_\xi(z)$ of any origin-symmetric star body $K$ with $C^\infty$ boundary. For $q \in \mathbb{C}$ and $\xi \in \mathbb{R}^n \setminus \{0\}$, we define

(6)
$$A_\xi^{(q)}(0) = \left\langle \frac{t_+^{-q-1}}{\Gamma(-q)}, A_\xi(t)\right\rangle.$$

If $m \in \mathbb{N}$ and $\operatorname{Re} q < m$, $q \ne 0, 1, 2, \ldots, m-1$, then

(7)
$$A_\xi^{(q)}(0) = \frac{1}{\Gamma(-q)} \int_0^1 t^{-q-1}$$
$$\cdot \left(A_\xi(t) - A_\xi(0) - tA_\xi'(0) - \cdots - \frac{t^{m-1}}{(m-1)!}A_\xi^{(m-1)}(0)\right) dt$$
$$+ \frac{1}{\Gamma(-q)} \int_1^\infty t^{-q-1} A_\xi(t)\, dt + \frac{1}{\Gamma(-q)} \sum_{k=0}^{m-1} \frac{A^{(k)}(0)}{k!(k-q)},$$

and if $m - 1 < \operatorname{Re} q < m$, then

$$A_\xi^{(q)}(0) = \frac{1}{\Gamma(-q)} \int_0^\infty t^{-q-1}$$
$$\cdot \left( A_\xi(t) - A_\xi(0) - tA_\xi'(0) - \cdots - \frac{t^{m-1}}{(m-1)!} A_\xi^{(m-1)}(0) \right) dt.$$

Furthermore, we deduce that $q \mapsto A_\xi^{(q)}(0)$ is an analytic function on $\mathbb{C}$ with

$$A_\xi^{(k)}(0) = (-1)^k \left. \frac{\partial^k}{\partial t^k} A_\xi(t) \right|_{t=0} \quad \text{for } k = 0, 1, 2, \ldots.$$

Since $K$ is symmetric, the function $t \mapsto A_\xi(t)$ is even, and for every even $m$,

$$(8) \quad A_\xi^{(q)}(0) = \frac{1}{\Gamma(-q)} \int_0^\infty t^{-q-1} \left( A_\xi(t) - \sum_{j=0}^{(m-2)/2} \frac{t^{2j}}{(2j)!} A_\xi^{(2j)}(0) \right) dt,$$

whenever $m - 2 < \operatorname{Re} q < m$.

We remark that (7) (and so also (8)) was deduced from (4) and (5), and the fact that $A_\xi^{(q)}(0)$ is an analytic function in $q$ was deduced from the fact that $\langle t_+^\lambda / \Gamma(\lambda + 1), \cdot \rangle$ is an analytic distribution. We could also have used (7) as the definition of $A_\xi^{(q)}(0)$, $q \neq 1, 2, \ldots$, from which it is easily seen that $A_\xi^{(q)}(0)$ can be extended analytically to all of $\mathbb{C}$.

## 4. Proof of Theorems 1 and 2

Theorem 1 will be an easy consequence of Theorem 2. To prove the latter, we shall need the following lemma from [K1]. For the sake of completeness we include a proof.

LEMMA 5. *For every even test function $\varphi \in \mathcal{S}$, $\xi \in S^{n-1}$, and $-1 < q < 0$,*

$$\int_{\mathbb{R}^n} |\langle \xi, x \rangle|^{-q-1} \varphi(x) \, dx = \frac{-1}{2\Gamma(1+q) \sin \frac{q\pi}{2}} \int_{-\infty}^\infty |t|^q \hat{\varphi}(t\xi) \, dt.$$

*Proof.* Using (1) and (3), we deduce that

$$\int_{\mathbb{R}^n} |\langle \xi, x \rangle|^{-q-1} \varphi(x) \, dx = \int_\mathbb{R} |t|^{-q-1} \int_{\langle x, \xi \rangle = t} \varphi(x) \, dx \, dt$$
$$= \langle |t|^{-q-1}, \mathcal{R}\varphi(\xi; t) \rangle$$
$$= \left\langle \frac{-1}{2\Gamma(1+q) \sin \frac{q\pi}{2}} (|s|^q)^\wedge(t), \mathcal{R}\varphi(\xi; t) \right\rangle$$
$$= \frac{-1}{2\Gamma(1+q) \sin \frac{q\pi}{2}} \langle |s|^q, \hat{\varphi}(s\xi) \rangle. \qquad \square$$

*Proof of Theorem* 2. Suppose that $-1 < q < 0$. The function $A_\xi(z) = \int_{\langle x,\xi \rangle = z} \chi(\|x\|)\, dx$ is even. Applying Fubini's theorem and passing to spherical coordinates, we get

$$A_\xi^{(q)}(0) = \frac{1}{2\Gamma(-q)} \int_{-\infty}^{\infty} |z|^{-q-1} A_\xi(z)\, dz$$

$$= \frac{1}{2\Gamma(-q)} \int_{\mathbb{R}^n} |\langle x, \xi \rangle|^{-q-1} \chi(\|x\|)\, dx$$

$$= \frac{1}{2\Gamma(-q)} \int_{S^{n-1}} |\langle \theta, \xi \rangle|^{-q-1} \int_0^\infty r^{n-q-2} \chi(r\|\theta\|)\, dr\, d\theta$$

$$= \frac{1}{2(n-q-1)\Gamma(-q)} \int_{S^{n-1}} |\langle \theta, \xi \rangle|^{-q-1} \|\theta\|^{-n+q+1}\, d\theta.$$

We now consider $A_\xi^{(q)}(0)$ as a function of $\xi \in \mathbb{R}^n \setminus \{0\}$. By Lemma 5, for every even test function $\varphi \in \mathcal{S}$,

(9)
$$\langle A_\xi^{(q)}(0), \varphi(\xi) \rangle = \frac{1}{2(n-q-1)\Gamma(-q)} \int_{S^{n-1}} \|\theta\|^{-n+q+1} \int_{\mathbb{R}^n} |\langle \theta, \xi \rangle|^{-q-1} \varphi(\xi)\, d\xi$$

$$= \frac{-1}{4(n-q-1)\Gamma(-q)\Gamma(q+1)\sin\frac{q\pi}{2}} \int_{S^{n-1}} \|\theta\|^{-n+q+1} \int_{-\infty}^\infty |t|^q \widehat{\varphi}(t\theta)\, dt\, d\theta$$

$$= \frac{\cos\frac{q\pi}{2}}{\pi(n-q-1)} \langle (\|x\|^{-n+q+1})^\wedge(\xi), \varphi(\xi) \rangle,$$

where the last equation follows from the property of the gamma function that $\Gamma(-q)\Gamma(q+1) = -\pi/\sin(q\pi)$ and a simple calculation

$$\langle (\|x\|^{-n+q+1})^\wedge(\xi), \varphi(\xi) \rangle = \int_{\mathbb{R}^n} \|x\|^{-n+q+1} \widehat{\varphi}(x)\, dx$$

$$= \int_{S^{n-1}} \|\theta\|^{-n+q+1} \int_0^\infty t^q \widehat{\varphi}(t\theta)\, dt\, d\theta$$

(note that the function $\|x\|^{-n+q+1}$ is locally integrable on $\mathbb{R}^n$ because $-1 < q < 0$).

Since (9) holds for every even test function $\varphi$, Theorem 2 is proved when $-1 < q < 0$.

In order to prove the theorem for other values of $q$, we first observe that $(\|x\|^{-n+q+1})^\wedge$ is an analytic distribution (with respect to $q$) on $\{q \in \mathbb{C} : \operatorname{Re} q > -1\}$. It follows that for every even test function $\varphi \in \mathcal{S}$, the functions $q \mapsto \langle A_\xi^q(0), \varphi \rangle$ (see §3) and

$$q \mapsto \left\langle \frac{\cos\frac{q\pi}{2}}{\pi(n-q-1)} (\|x\|^{-n+q+1})^\wedge(\xi), \varphi \right\rangle$$

are analytic on the connected region $\{q \in \mathbb{C}: \operatorname{Re} q > -1, q \neq n - 1\}$. These functions coincide on the interval $-1 < q < 0$, so they coincide on $\{q \in \mathbb{C}: \operatorname{Re} q > -1, q \neq n - 1\}$. Since $\varphi$ is an arbitrary even test function, we have proved Theorem 2. $\square$

*Proof of Theorem* 1. If $k$ is even, the theorem follows immediately from Theorem 2 (with $q = k$) and the equation $\sin \frac{(k+1)\pi}{2} = (-1)^{k/2}$.

If $k$ is odd, both sides of the equation in Theorem 2 vanish. Assuming that $q \neq 0, 1, 2, \ldots$ in this equation, multiply both sides of the equation by $\Gamma(-q)$, and then take the limit as $q \to k$.

By (8) with $m = k + 1$,

$$\lim_{q \to k} \Gamma(-q) A_\xi^{(q)}(0) = \int_0^\infty t^{-k-1} \left( A_\xi(t) - \sum_{j=0}^{(k-1)/2} \frac{t^{2j}}{(2j)!} A_\xi^{(2j)}(0) \right) dt.$$

We also observe, using $\Gamma(\lambda + 1) = \lambda \Gamma(\lambda)$, that

$$\lim_{q \to k} \Gamma(-q) \sin \frac{(q+1)\pi}{2} = \lim_{q \to k} \Gamma(-q) \sin \left( \frac{(q+1)\pi}{2} - \frac{(k+1)\pi}{2} \right) (-1)^{(k+1)/2}$$

$$= \lim_{q \to k} \frac{\Gamma(-q + k + 1)}{(-q)(1-q)\cdots(k-q)} \sin \frac{(q-k)\pi}{2} (-1)^{(k+1)/2}$$

$$= -\frac{\pi}{2} (-1)^{(k+1)/2} (-1)^k \frac{1}{k!} = \frac{\pi}{2} (-1)^{(k+1)/2} \frac{1}{k!}.$$

The statement of Theorem 1(b) follows. $\square$

WESTERN WASHINGTON UNIVERSITY, BELLINGHAM, WA
*E-mail address*: gardner@baker.math.wwu.edu

UNIVERSITY OF TEXAS AT SAN ANTONIO, SAN ANTONIO, TX
*E-mail address*: koldobsk@sphere.math.utsa.edu

TEXAS A&M UNIVERSITY, COLLEGE STATION, TX
*E-mail address*: schlump@math.tamu.edu